\documentclass[aos,MSNbibl,nameyear,dvips]{arximspdf}
\usepackage{graphicx}
\doi{10.1214/14-AOS1266} 
\volume{43}
\issue{1}
\pubyear{2015}
\firstpage{90}
\lastpage{101}
\docsubty{FLA}

\makeatletter
\def\textsc{\mathrm}
\newproclaim{remark}{Remark}

\newtheorem{theorem}{Theorem}

\newcommand{\Ga}{\Gamma}
\newcommand{\de}{\delta}
\newcommand{\ka}{\kappa}
\newcommand{\ep}{\varepsilon}
\newcommand{\La}{\Lambda}
\makeatother

\begin{document}
\begin{frontmatter}

\title{A new permutation test statistic for complete block designs}
\runtitle{Permutation tests for complete block design}

\begin{aug}
\author[A]{\fnms{Inga}~\snm{Samonenko}\ead[label=e1]{isamonenko@yahoo.com}\thanksref{E1}}
\and
\author[A]{\fnms{John}~\snm{Robinson}\corref{}\ead[label=e2]{john.robinson@sydney.edu.au}\thanksref{E2}}
\runauthor{I. Samonenko and J. Robinson}
\thankstext{E1}{Supported by UPA Scholarship.}
\thankstext{E2}{Supported by ARC DP0773345.}
\affiliation{University of Sydney}
\address[A]{School of Mathematics and Statistics\\
University of Sydney\\
NSW 2006\\
Australia\\
\printead{e1}\\
\phantom{E-mail:\ }\printead*{e2}}

\end{aug}

\received{\smonth{12} \syear{2013}}
\revised{\smonth{4} \syear{2014}}

%
\begin{abstract}
We introduce a nonparametric test statistic for the permutation test
in complete block designs. We find the region in which the statistic
exists and consider particularly its properties on the boundary of the
region. Further, we prove that saddlepoint approximations for tail
probabilities can be obtained inside the interior of this region.
Finally, numerical examples are given showing that both accuracy and
power of the new statistic improves on these properties of the
classical $F$-statistic under some non-Gaussian models and equals them
for the Gaussian case.
\end{abstract}

%
\end{frontmatter}
%
\section{Introduction}\label{sec1}
Randomized designs and permutation tests originated in the work of
\citet{Fish}. \citet{KR11} obtained theorems on the distribution of a
general likelihood ratio like statistic under weak conditions and
applied these to the one-way or $k$-sample permutation tests, obtaining
saddlepoint approximations generalizing the Lugananni--Rice and
Barndorff--Nielsen approximations for one-dimensional means. Here, we
use their general result and apply their approach to permutation tests
for complete block designs, paying particular attention to the region
in which the statistic exists and in the interior of which saddlepoint
approximations can be obtained. This interior is the admissible domain,
following \citet{BR65}. We examine the properties of the test statistic
in this region and on its boundary, and obtain results on the relative
errors of saddlepoint approximations inside the admissible domain. We
also give numerical results for comparisons of the new statistic with
the commonly used $F$-statistic which demonstrate the accuracy of the
saddlepoint approximation and show, for long tailed error
distributions, an improvement in power relative to the $F$-statistic
with no loss of power for near normal errors.

A randomized complete block design is used to compare the effect of $k$
different treatments in $b$ blocks, usually selected to reduce the
variation within subunits of the block. The analysis of variance is
used to test the null hypothesis that the treatments have the same
effect, with the test statistic $F$, the ratio of the treatment and
error mean squares.
Under the assumption that the errors are normally distributed, the null
distribution of $F$ is the $\mathit{F}$ distribution with $k-1$ and
$(k-1)(b-1)$-degrees of freedom and the $F$ test is equivalent to an
unconditional likelihood ratio test.

The random assignment of $k$ treatments to each block allows us to use
a permutation test based on means which is distribution-free and does
not rely either on the assumption of normality or on asymptotics. This
test can be performed using the $F$-statistic and will have correct
size, conditionally on the order statistics in each block, and so
unconditionally, for any distribution of errors under the null
hypothesis of no treatment effects in a standard two-way model or for a
model based on randomization prior to the experiment. Under the null
hypothesis, the permutation distribution of this statistic can be
calculated exactly by evaluating all possible values of the test
statistic under permutations in each block and taking these as
equi-probable. When this is numerically infeasible, Monte Carlo methods
are widely used to approximate the exact distribution by using a large
random sample of the possible permutations. A chi-squared distribution
with $(k-1)$-degrees of freedom or an $F$ distribution with $(k-1)$ and
$(k-1)(b-1)$-degrees of freedom are asymptotic approximations to the
distribution of the permutation test statistic under mild conditions on
moments. If the observations are not normally distributed and if the
number of blocks is not large, then the central limit theorem will not
guarantee a good approximation and the test will not have the
optimality properties that might be expected under normality.

We propose a likelihood ratio like statistic in place of $\mathit{F}$,
based on exponential tilting. We show that this statistic can be
calculated on the admissible domain, an open convex set, the closure of
which contains the support of the treatment means. We consider the
boundary of the admissible domain and show that the statistic can be
obtained on the boundary as a limit which can be calculated using lower
dimensional versions of the statistic on lower dimensional versions of
the admissible domain. We then obtain saddlepoint approximations for
the tail probability of this statistic with relative errors of order
$1/n$ in the admissible domain, based on Theorems of \citet{KR11}. The
results generalize the saddlepoint approximations of \citet{R82} in the
case of permutation tests of paired units, which can be regarded as a
block design with blocks of size 2, where the admissible domain is the
interval between the mean of the absolute values of differences of the
pairs and its negative.

In the next section, we introduce the notation for a complete block
design, obtain the likelihood ratio like statistic and define its
admissible domain. In Section~\ref{sec3}, we describe the admissible domain
and give three theorems giving explicit results for the test statistic
on the boundary of the domain, with proofs given in Section~\ref{sec6}. In
Section~\ref{sec4}, we use the theorems of \citet{KR11} to show that tail
probabilities for the statistic under permutations can be approximated
in the admissible domain by an integral of a formal saddlepoint density
given in forms like those of Lugananni--Rice and Barndorff--Nielsen in
the one-dimensional case. In Section~\ref{sec5}, we present numerical
calculations illustrating the accuracy of the approximations compared
to those obtained using the standard test statistics and give power
comparisons showing an improvement in power over the standard $F$-test
for observations from long tailed distributions. The code used is
available from \url{http://www.maths.usyd.edu.au/u/johnr/BlockDesfns.R}.

\section{The test statistic \texorpdfstring{$\Lambda$}{Lambda} and its admissible domain}\label{sec2}

Let ($X_{ij}$) be a matrix of observed experimental values, normalized
to have row means zero, where $i = 1,\ldots, b$ is the block number and
$j=1, \ldots, k$ is the treatment number. Let matrix $A=(a_{ij})$ have
rows of the matrix ($X_{ij}$) each set in ascending order and let $A_i$
be its $i$th row.
Define the means $\bar{X}_j= \sum_{i=1}^b X_{ij}/b$ for $j=1,\ldots,k$
and let $\bar X=(\bar X_1,\ldots,\bar X_{k-1})^T$. Then, given $A$,
under the null hypothesis of equal treatment effects, the conditional
cumulative generating function for treatment means is
\[
b\ka(\tau)=\log\mathbf{E} \bigl(e^{\sum_{j=1}^k\tau_j\bar{X}_j} | A\bigr)=\sum
_{i=1}^b \log\mathbf{E} \bigl( e^{\sum_{j=1}^k\tau_j X_{ij}/b} |
A_i \bigr).
\]
Set $t_i = (\tau_i - \tau_k)/b$ for $i=1,\ldots,k-1$. Then we can reduce
the problem of defining the average cumulative generating function to a
$(k-1)$-dimensional one and write
%
\begin{equation}
\label{emcum} \ka(t)=\frac{1}{b}\sum_{i=1}^b
\log\frac{1}{k!}\sum_{\pi\in\Pi} e^{t^T a_{i\pi}},
\end{equation}
where $t=(t_1,\ldots,t_{k-1})^T$, $\Pi$ is the set of possible vectors
$(\pi(1), \ldots, \pi(k-1))$ obtained from the first $k-1$ elements of
all permutations of indices $\{1,\ldots,k\}$ and $a_{i\pi}=(a_{i\pi
(1)},\ldots, a_{i\pi(k-1)})^T$.

Consider the test statistic $\Lambda(\bar X)$, where
%
\begin{equation}
\label{emlam} \La(x) = \sup_{t} \bigl\{t^Tx-
\ka(t)\bigr\},
\end{equation}
for $t, x \in\mathbf{R}^{k-1}$. Let us define an admissible domain
$\Omega\subset\mathbf{R}^{k-1}$ as the set of all $x$ for which $t^T
x - \ka(t)$ attains its maximum. Then there exists a unique value $t_x$
such that
%
\begin{equation}
\label{emspe} \La(x) = t_x^Tx-\ka(t_x)\quad
\mbox{if and only if}\quad \ka'(t_x)=x,
\end{equation}
since $t^Tx-\ka(t)$ is strictly convex by noting that $ - \ka''(t)$ is
negative definite.

In the case $k=2$, the admissible domain is $(-\sum_{i=1}^b|a_{i1}-a_{i2}|,\sum_{i=1}^b|a_{i1}-a_{i2}|)$ and the
properties of $\Lambda$ and the saddlepoint approximation are discussed
for the two special cases of the binomial and the Wilcoxon signed-rank
test in \citet{JR}. The situation is more complex for $k>2$ and results
are given in the next section.

\begin{remark*}
Exact randomization tests have restricted
application to designed experiments. The only two designs for which we
know how to obtain a statistic of our form are the complete block
design considered here and the one-way or $k$-sample design considered
in \citet{KR11}. An extension to some other cases such as balanced
incomplete block designs or in testing for main effects using
restricted randomization as suggested by \citet{BM} may be possible but
do not seem to be straightforward.
\end{remark*}

\section{The properties of \texorpdfstring{$\Lambda$}{Lambda}}\label{sec3}

First, we will describe the admissible domain and give some results
which make it possible to calculate $\Lambda(x)$ on the boundary of the
domain where the solution of the saddlepoint equations (\ref{emspe})
does not exist. Let $\bar{A} = (\bar{A}_1,\ldots,\bar{A}_k)^T$ be a
vector of column means of $A$ and write $\bar{A}_\pi= (\bar{A}_{\pi
(1)}, \ldots, \bar{A}_{\pi(k-1)})^T$, for any $\pi\in\Pi$. Then the
support of $\bar X$ contains $\bar A_\pi$ and the set of $\bar A_\pi$,
for all $\pi\in\Pi$, is the set of extreme points of the convex hull of
the support of $\bar X$, which is a $(k-1)$-polytope $\mathcal{P}$.

\begin{theorem}\label{theo2}
The set $\Omega$ is the interior of the $(k-1)$-polytope $\mathcal{P}$.
\end{theorem}

%
\begin{theorem}\label{Th2}
The function $\La(x)$ is finite on the boundary of $\Omega$ and takes
its maximum value $\log k!$ at its extreme points.
\end{theorem}

The boundary of $\mathcal{P}$ consists of all $x \in\mathcal{P}$ for
which there exists an integer $l$ and distinct integers $s_1, \ldots,
s_{k-1}$ from the set $\{1, \ldots, k-1\}$ satisfying one of the equalities
%
\begin{equation}
\sum_{j=1}^{l} x_{s_j}=\sum
_{j=1}^{l} \bar{A}_j \quad\mbox{or}\quad
\sum_{j=1}^{l} x_{s_j}=\sum
_{j=1}^{l} \bar{A}_{k-j+1}.\label{bdyeqs}
\end{equation}

\begin{theorem}\label{Th3}
On the boundary of $\Omega$ corresponding to the value $l$ we have
\[
\La(x) = \La_1(x) + \La_2(x) + \log\pmatrix{k
\cr
l},
\]
for
\[
\La_1(x) = \sup_{u_1,\ldots,u_{l-1}} \Biggl( \sum
_{j=1}^{l-1} x_{s_j} u_j-
\frac{1}{b}\sum_{i=1}^b \log
\frac{1}{l!} \sum_{\hat{\pi}_1 \in\hat{\Pi}_1} e^{\sum_{j=1}^{l-1} a_{i\hat{\pi}_1(j)} u_j} \Biggr)
\]
and
\[
\La_2(x) = \sup_{u_{l+1},\ldots,u_{k-1}} \Biggl( \sum
_{j=l+1}^{k-1} x_{s_j} u_j -
\frac{1}{b}\sum_{i=1}^b \log
\frac{1}{(k-l)!}\sum_{\hat{\pi
}_2 \in\hat{\Pi}_2}e^{\sum_{j=l+1}^{k-1} a_{i\hat{\pi}_2(j)} u_j} \Biggr),
\]
where $\hat{\Pi}_1$ and $\hat{\Pi}_2$ are sets of all permutations
$\hat
{\pi}_1$ and $\hat{\pi}_2$ of integers $\{ 1, \ldots, l\}$ and $\{ l+1,
\ldots, k\}$, respectively.
\end{theorem}

\begin{figure}

\includegraphics{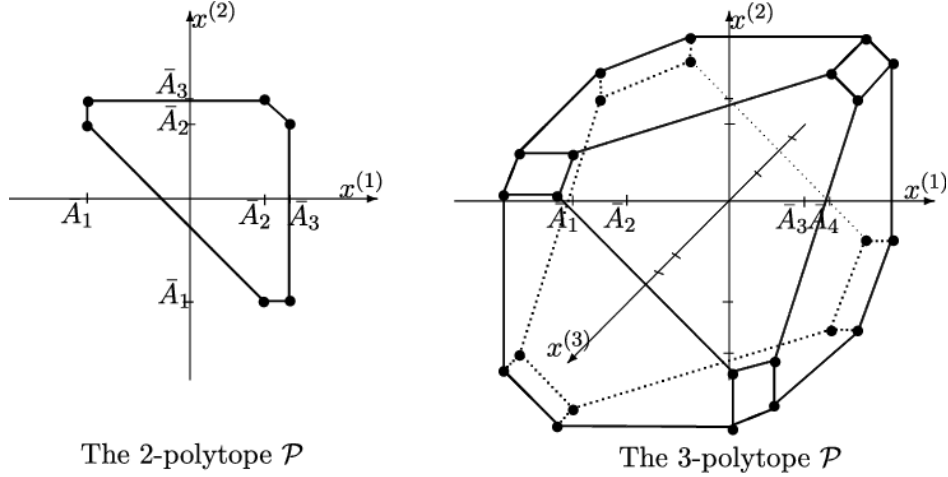}

  \caption{Examples of the admissible domain in the
  cases $k=3$ and $k=4$.}\label{fig1}
\end{figure}

\begin{remark*} The result of Theorem~\ref{Th3} demonstrates that the
boundary of
$\Omega\subset\mathbf{R}^{k-1}$ consists of lower dimensional
polytopes, each made up of a cross product of two sets of dimension
$l-1$ and $k-l-1$, for $l=1,\ldots,k-1$. These correspond to the
restriction of the permutations in each block to the smallest or
largest $l-1$ elements of the block and their complements. The
functions $\Lambda_1$ and $\Lambda_2$ are defined on these subsets as
is $\Lambda$ in (\ref{emlam}). To illustrate this, in Figure \ref{fig1}
we have given two
diagrams showing the polytope $\mathcal{P}$ for the cases $k=3$ and
$k=4$. In the first picture, we have 6 vertexes and 6 sides with
boundaries made up of lines representing the dimension reduction to one
dimension. In this case, one of $\Lambda_1$ and $\Lambda_2$ is
identically zero. In the second picture, the two-dimensional boundaries
are either six-sided, corresponding to one of $\Lambda_1$ and $\Lambda
_2$ being identically zero, and the other a two-dimensional function,
or are rectangles corresponding to both $\Lambda_1$ and $\Lambda_2$
being one-dimensional functions.
\end{remark*}

\section{Saddlepoint approximations for \texorpdfstring{$\Lambda$}{Lambda}}\label{sec4}
Consider $\textsc{P}(\La(\bar{X}) \ge u^2/2)$, where $\textsc{P}$
denotes the conditional distribution given $A$, and define
\[
r(x) = e^{-b \La(x)} (2\pi/ b)^{-k/2}\bigl|\ka''(t_x)\bigr|^{-1/2}
\]
for $x\in\Omega$. In the case of identically distributed random
vectors ${X}_i=(X_{i1},\break  \ldots,  X_{i,k-1}) \in\Omega$, $i=1,\ldots,N$
with known densities, $r(\bar{X})$ is a saddlepoint density
approximation for $\bar{X}$, obtained by \citet{BR65}. In our case, the
lack of a density requires the application of Theorem~1 of \citet{KR11}.
We consider
\[
\textsc{P}\bigl(\La(\bar{X}) \ge u^2/2\bigr)=\textsc{P}(\bar X\in
\mathcal{F}),
\]
where $\mathcal{F} =\{ x\dvtx \La(x) \ge u^2/2 \}$ and $u\in\Omega
_{-\ep}
= (\Omega_{\ep}^c)^c = \{x \in\mathbf{R}^{k-1}, y \in\Omega^c\dvtx |x-y|<\ep\}^c$. Whenever $\mathcal{F}^c \subset\Omega_{-\ep}$, our
case must only meet the necessary conditions (A1)--(A4) stated in
\citet{KR11}. The cumulative generating function (\ref{emcum}) exists
throughout $\mathbf{R}^{k-1}$, therefore, the first condition is met.
The average variance $\ka''(t)$ is a positive definite matrix which
equals the identity matrix at the origin. Thus, the second condition is
met. The third condition only requires the existence of some moments
and the fourth is a smoothness condition, which we assume holds. It
will hold, for example, if the observations are from a distribution
with a continuous component. Thus we can apply Theorems 1 and 2 of
\citet{KR11} as in that paper to get the following result.

\begin{theorem}\label{Th4}
For $\varepsilon>0$ and $u^2/2<\log k-\varepsilon$, under conditions
\textup{(A1)--(A4)} of \citet{KR11},
%
\begin{eqnarray}\qquad
\mathbf{P}\bigl(\Lambda(\bar{X}) \ge u^2/2\bigr) &=&
Q_{k-1}\bigl(b u^{2}\bigr) \bigl(1+O(1/b)\bigr)+
\frac{c_b}{b}u^{k-1}e^{-{bu^2}/{2}}\frac{G(u)-1}{u^2},\label{LR}
\\
\mathbf{P}\bigl(\Lambda(\bar{X}) \ge u^2/2\bigr) &=&
Q_{k-1}\bigl(b u^{*2}\bigr) \bigl(1+O(1/b)
\bigr),\label{BN}
\end{eqnarray}
where
\[
Q_{k-1} (x) = \textsc{P}\bigl(\mathcal{X}^2_{k-1}
\ge x\bigr)=\frac
{1}{2^{(k-1)/2} \Ga((k-1)/2)} \int_x^{\infty}z^{(k-1)/2-1}e^{-z/2}
\,dz,
\]
$ u^*= u - \log(G(u))/bu$, $ c_b = b^{{(k-1)}/{2}}/2^{{(k-3)}/{2}}
\Ga
(\frac{k-1}{2})$,
\[
\de(u,s) = \frac{\Ga({(k-1)}/{2}) |\ka''(t_x)|^{-{1}/{2}}
|\ka''(0)| ^{{1}/{2}} r^{k-2}} {2\pi^{{(k-1)}/{2}} u^{k-3}
|s^T\ka''(0)^{{1}/{2}}t_x|}
\]
and
\[
G(u) = \int_{S_{k-1}} \de(u,s) \,ds,
\]
for $S_{k-1}$ the $k-1$-dimensional unit sphere centered at zero, and
where, for each $s\in S_{k-1}$, $r$ is chosen so $\La(r\ka
''(0)^{1/2}s) = u^2/2$. Here, $t_x$ is a solution to (\ref{emspe}) at
the point $x=r \ka''(0)^{1/2}s$.
\end{theorem}

We note that the constraint $u^2/2<\log k-\varepsilon$, ensures that the
level set of $\Lambda(x)$ corresponding to $u$ lies entirely in
$\Omega
$, since the minimum value of $\Lambda(x)$ for an $x$ on the boundary
occurs for $l=1$ and $\Lambda_1(x)=\Lambda_2(x)=0$ in Theorem~\ref{Th3}. The
remainder of the proof then follows in the same way as in Theorem~2 of
\citet{KR11}, so we omit it.
The theorem gives approximations of the tail probabilities of the test
statistic $\La$ under permutations, in forms like those of
Lugananni--Rice and Barndorff--Nielsen in the one-dimensional case.

\begin{table}
\caption{Accuracy for exponentially squared errors, $b=10$ and $k=4$}
\label{T2}
\begin{tabular*}{\textwidth}{@{\extracolsep{\fill}}lccccc@{}}
\hline
\multicolumn{1}{@{}l}{$\bolds{u}$} &
\multicolumn{1}{c}{\textbf{0.6}} &
\multicolumn{1}{c}{\textbf{0.8}} &
\multicolumn{1}{c}{\textbf{1.0}} &
\multicolumn{1}{c}{\textbf{1.2}} &
\multicolumn{1}{c@{}}{\textbf{1.4}} \\
\hline
\multicolumn{1}{@{}l}{MC $F$} &
\multicolumn{1}{c}{0.3353} &
\multicolumn{1}{c}{0.1042} &
\multicolumn{1}{c}{0.0257} &
\multicolumn{1}{c}{0.0047} &
\multicolumn{1}{c@{}}{0.0003} \\
\multicolumn{1}{@{}l}{$F$} &
\multicolumn{1}{c}{0.3286} &
\multicolumn{1}{c}{0.1193} &
\multicolumn{1}{c}{0.0342} &
\multicolumn{1}{c}{0.0083} &
\multicolumn{1}{c@{}}{0.0018} \\
\multicolumn{1}{@{}l}{MC $\Lambda$} &
\multicolumn{1}{c}{0.4494} &
\multicolumn{1}{c}{0.1777} &
\multicolumn{1}{c}{0.0455} &
\multicolumn{1}{c}{0.0071} &
\multicolumn{1}{c@{}}{0.0008} \\
\multicolumn{1}{@{}l}{SP LR} &
\multicolumn{1}{c}{0.4179} &
\multicolumn{1}{c}{0.1645} &
\multicolumn{1}{c}{0.0480} &
\multicolumn{1}{c}{0.0070} &
\multicolumn{1}{c@{}}{0.0007} \\
\multicolumn{1}{@{}l}{SP BN} &
\multicolumn{1}{c}{0.4059} &
\multicolumn{1}{c}{0.1565} &
\multicolumn{1}{c}{0.0441} &
\multicolumn{1}{c}{0.0066} &
\multicolumn{1}{c@{}}{0.0007} \\
\hline
\end{tabular*}\vspace*{9pt}
\end{table}

\begin{table}[b]\vspace*{9pt}
\caption{Accuracy for exponential squared errors, $b=5$ and $k=3$}
\label{T3}
\begin{tabular*}{\textwidth}{@{\extracolsep{\fill}}lccccc@{}}
\hline
\multicolumn{1}{@{}l}{$\bolds{u}$} &
\multicolumn{1}{c}{\textbf{0.6}} &
\multicolumn{1}{c}{\textbf{0.8}} &
\multicolumn{1}{c}{\textbf{1.0}} &
\multicolumn{1}{c}{\textbf{1.2}} &
\multicolumn{1}{c@{}}{\textbf{1.4}} \\
\hline
\multicolumn{1}{@{}l}{MC $F$} &
\multicolumn{1}{c}{0.4620} &
\multicolumn{1}{c}{0.2448} &
\multicolumn{1}{c}{0.1277} &
\multicolumn{1}{c}{0.0653} &
\multicolumn{1}{c@{}}{0.0427} \\
\multicolumn{1}{@{}l}{$F$} &
\multicolumn{1}{c}{0.4441} &
\multicolumn{1}{c}{0.2603} &
\multicolumn{1}{c}{0.1434} &
\multicolumn{1}{c}{0.0767} &
\multicolumn{1}{c@{}}{0.0408} \\
\multicolumn{1}{@{}l}{MC $\Lambda$} &
\multicolumn{1}{c}{0.5118} &
\multicolumn{1}{c}{0.3133} &
\multicolumn{1}{c}{0.1578} &
\multicolumn{1}{c}{0.0774} &
\multicolumn{1}{c@{}}{0.0270} \\
\multicolumn{1}{@{}l}{SP LR} &
\multicolumn{1}{c}{0.5011} &
\multicolumn{1}{c}{0.2988} &
\multicolumn{1}{c}{0.1563} &
\multicolumn{1}{c}{0.0736} &
\multicolumn{1}{c@{}}{0.0278} \\
\multicolumn{1}{@{}l}{SP BN} &
\multicolumn{1}{c}{0.4950} &
\multicolumn{1}{c}{0.2917} &
\multicolumn{1}{c}{0.1500} &
\multicolumn{1}{c}{0.0687} &
\multicolumn{1}{c@{}}{0.0255} \\
\hline
\end{tabular*}
\end{table}

\section{Numerical results}\label{sec5}

\subsection{Accuracy}\label{sec5.1}
For each of the simulation experiments, we obtained a single matrix $A$
by sampling from a distribution, that of squared exponential random
variables for our Tables~\ref{T2}, \ref{T3} and \ref{T4}. Then we used
100,000 replicates of random permutations of each block to obtain Monte
Carlo approximations to the tail probabilities of the permutation tests
for the statistics $F$ and $\Lambda$, shown as $\mathrm{MC}$ $F$ and $\mathrm{MC}$
$\Lambda$ in the tables. We compared these to the tail probabilities
from the $F$ distribution for the $F$-statistic and to the saddlepoint
approximations for the $\Lambda$ statistic obtained using formulas
(\ref
{LR}) (SP LR) and (\ref{BN}) (SP BN), respectively, with $100$ Monte
Carlo samples used to approximate integrals on the sphere $S_{k-1}$, as
in the Remark in Section~2 of \citet{KR11}. We also used the method from
\citet{Genz}, for approximation of the integral on the sphere, obtaining
effectively the same accuracy as with Monte Carlo sampling.

From Tables~\ref{T2} and \ref{T3}, we note that the accuracy is high
for the $\Lambda$ test, even for only 5 blocks of size 3. We note that
for Table~\ref{T3} the theorem holds for $u$ less than $\sqrt{2\log
3}=1.48$, so we are restricted to this region. Results from other
simulations show even greater accuracy under normal errors or errors
that are not from long tailed distributions.

The $F$-statistic has less accuracy in the tails, partly because the
$F$-statistic approximates the average of tail probabilities
conditioned on the matrix $A$, using $100\mbox{,}000$ permutations for each
$A$, so that even in the case of normal errors, it may not agree with
the conditional tail probabilities approximated by $\mathrm{MC}$ $F$-values
from the tables of this section, which give proportions in the tails
obtained from $100\mbox{,}000$ Monte Carlo simulations from a particular
sample and is an approximation of the conditional distribution. To
consider the accuracy of the unconditional test, we obtained 1000
replicates from each of a normal and exponential squared distribution,
obtained tail probabilities for these from the permutation test for the
$F$-statistic, averaged these over the 1000 replicates and compared
these approximations to the $F$ distribution. In the normal case, the
results were very accurate, essentially replicating the theoretical
results, as expected, and for the squared exponential case the results
are given in Table~\ref{T4} indicating that errors remain
unsatisfactory in the tails.

\begin{table}
\caption{Average of 1000 permutation test results for exponential
squared errors with $b=10$ and $k=4$ compared to the $F$ distribution}
\label{T4}
\begin{tabular*}{\textwidth}{@{\extracolsep{\fill}}lccccc@{}}
\hline
\multicolumn{1}{@{}l}{$\bolds{u}$} &
\multicolumn{1}{c}{\textbf{0.6}} &
\multicolumn{1}{c}{\textbf{0.8}} &
\multicolumn{1}{c}{\textbf{1.0}} &
\multicolumn{1}{c}{\textbf{1.2}} &
\multicolumn{1}{c@{}}{\textbf{1.4}} \\
\hline
\multicolumn{1}{@{}l}{$E$ MC $F$} &
\multicolumn{1}{c}{0.3175} &
\multicolumn{1}{c}{0.0836} &
\multicolumn{1}{c}{0.0171} &
\multicolumn{1}{c}{0.0032} &
\multicolumn{1}{c@{}}{0.0005} \\
\multicolumn{1}{@{}l}{$F$} &
\multicolumn{1}{c}{0.3286} &
\multicolumn{1}{c}{0.1193} &
\multicolumn{1}{c}{0.0342} &
\multicolumn{1}{c}{0.0083} &
\multicolumn{1}{c@{}}{0.0018} \\
\hline
\end{tabular*}
\end{table}

\subsection{Power results for the saddlepoint approximations}\label{sec5.2}
We compare the $F$-statistic and the saddle point approximations using
(\ref{LR}) and (\ref{BN}) using $100$ Monte Carlo uniform samples from
$S_k$. There were $2000$ samples with errors drawn from the exponential
and the exponential squared distributions, and for each of these
$p$-values were calculated using the saddlepoint approximations for the
$\Lambda$ statistic obtained using formulas (\ref{LR}) (PowerLR) and
(\ref{BN}) (PowerBN), respectively, and using $10\mbox{,}000$ permutations to
approximate the $p$-values for the $F$-statistic, for a design with $10$
blocks of size $4$. We selected treatment effects $\mu$, as $(0, 0, 0,
0), (-1/5, 0, 0, 1/5), \ldots, (-9/5, 0, 0, 9/5)$.

Under the exponential distribution, in Table~\ref{T5}, the
$\Lambda$-statistic gives a slight increase in power compared to
$F$-statistic
for small $\sum\mu^2$ and under the exponentially squared distribution,
in Table~\ref{T6}, the $\Lambda$-statistic gives a substantial increase
in power compared to $F$-statistic for moderate values of $\sum\mu^2$.
In both cases there is no difference for higher powers. We note that
the tests have essentially equal power up to computational accuracy
under the Normal, Uniform and Gamma (shape parameter 5) distributions.
The increase in power becomes noticeable in long tail distributions
like Exponential, Exponential Squared, Gamma (shape parameter 0.5)
distributions.

\begin{table}
\caption{Power calculation under the exponential distribution}
\label{T5}
\begin{tabular*}{\textwidth}{@{\extracolsep{\fill}}lcccccccc@{}}
\hline
\multicolumn{1}{@{}l}{$\boldsymbol{\sum}\bolds{\mu^2}$} & \multicolumn{1}{c}{\textbf{0.0}} &
\multicolumn{1}{c}{\textbf{0.04}} & \multicolumn{1}{c}{\textbf{0.16}} &
 \multicolumn
{1}{c}{\textbf{0.36}} & \multicolumn{1}{c}{\textbf{0.64}} &
\multicolumn{1}{c}{\textbf{1.0}} & \multicolumn{1}{c}{\textbf{1.44}} & \multicolumn
{1}{c@{}}{\textbf{1.96}}\\
\hline
\multicolumn{1}{@{}l}{PowerF} & \multicolumn{1}{c}{0.056} & \multicolumn
{1}{c}{0.081} & \multicolumn{1}{c}{0.209} & \multicolumn{1}{c}{0.424} &
\multicolumn{1}{c}{0.668} &
\multicolumn{1}{c}{0.861} & \multicolumn{1}{c}{0.982} & \multicolumn
{1}{c@{}}{1}\\
\multicolumn{1}{@{}l}{PowerLR} & \multicolumn{1}{c}{0.049} &
\multicolumn
{1}{c}{0.091} & \multicolumn{1}{c}{0.235} & \multicolumn{1}{c}{0.451} &
\multicolumn{1}{c}{0.679} &
\multicolumn{1}{c}{0.862} & \multicolumn{1}{c}{0.981} & \multicolumn
{1}{c@{}}{1}\\
\multicolumn{1}{@{}l}{PowerBN} & \multicolumn{1}{c}{0.053} &
\multicolumn
{1}{c}{0.093} & \multicolumn{1}{c}{0.238} & \multicolumn{1}{c}{0.455} &
\multicolumn{1}{c}{0.681} & \multicolumn{1}{c}{0.866} & \multicolumn
{1}{c}{0.982} & \multicolumn{1}{c@{}}{1} \\
\hline
\end{tabular*}
\end{table}

\begin{table}[b]
\caption{Power calculation under the exponentially squared distribution}
\label{T6}
\begin{tabular*}{\textwidth}{@{\extracolsep{\fill}}lcccccccc@{}}
\hline
\multicolumn{1}{@{}l}{$\bolds{\sum\mu^2}$} & \multicolumn{1}{c}{\textbf{0.0}} &
\multicolumn{1}{c}{\textbf{0.04}} & \multicolumn{1}{c}{\textbf{0.16}} &
 \multicolumn
{1}{c}{\textbf{0.36}} & \multicolumn{1}{c}{\textbf{0.64}} &
\multicolumn{1}{c}{\textbf{1.0}} & \multicolumn{1}{c}{\textbf{1.44}} & \multicolumn
{1}{c@{}}{\textbf{1.96}}\\
\hline
\multicolumn{1}{@{}l}{PowerF} & \multicolumn{1}{c}{0.051} & \multicolumn
{1}{c}{0.101} & \multicolumn{1}{c}{0.230} & \multicolumn{1}{c}{0.490} &
\multicolumn{1}{c}{0.711} &
\multicolumn{1}{c}{0.840} & \multicolumn{1}{c}{0.987} & \multicolumn
{1}{c@{}}{1}\\
\multicolumn{1}{@{}l}{PowerLR} & \multicolumn{1}{c}{0.057} &
\multicolumn
{1}{c}{0.169} & \multicolumn{1}{c}{0.319} & \multicolumn{1}{c}{0.545} &
\multicolumn{1}{c}{0.727} &
\multicolumn{1}{c}{0.832} & \multicolumn{1}{c}{0.976} & \multicolumn
{1}{c@{}}{1}\\
\multicolumn{1}{@{}l}{PowerBN} & \multicolumn{1}{c}{0.063} &
\multicolumn
{1}{c}{0.178} & \multicolumn{1}{c}{0.328} & \multicolumn{1}{c}{0.550} &
\multicolumn{1}{c}{0.731} &
\multicolumn{1}{c}{0.832} & \multicolumn{1}{c}{0.977} & \multicolumn
{1}{c@{}}{1}\\
\hline
\end{tabular*}
\end{table}

\section{Proofs of theorems of Section~\texorpdfstring{\protect\ref{sec3}}{3}}\label{sec6}
\mbox{}
\begin{pf*}{Proof of Theorem~\ref{theo2}}
From (\ref{emspe}), $\Omega=\{x\dvtx \kappa'(t)=x,\mbox{ for some }t\in
\mathbf{R}^{k-1}\}$, the image of $\kappa'(\cdot)$.
Using equation (\ref{emcum}) we get
the $j$th component of $\ka'(t)$,
\[
\kappa'_j(t)=\frac{1}{b}\sum
_{i=1}^b \frac{ \sum_{\pi\in\Pi}
a_{i\pi
(j)} \exp(t^T a_{i\pi})}{\sum_{\pi\in\Pi}\exp(t^T a_{i\pi})}.
\]
Here, $a_{i\pi(j)}$ is the $j$th component of $a_{i\pi}$ and
\[
\bar A_1= \frac{1}{b}\sum_{i=1}^b
a_{i1}<\frac{1}{b}\sum_{i=1}^b
\frac{
\sum_{\pi\in\Pi} a_{i\pi(j)} \exp(t^T a_{i\pi})}{\sum_{\pi\in
\Pi
}\exp(t^T a_{i\pi})} < \frac{1}{b}\sum_{i=1}^b
a_{ik} = \bar{A}_k.
\]
Using the same approach, we can conclude that for all distinct integers
$j_1$, $j_2, \ldots, j_{k-1}$ taking values $1, 2, \ldots, k-1$, and
for $l=1,\ldots,k-1$,
%
\begin{equation}\label{inteqs}
\sum_{j=1}^{l}\bar{A}_{j}<
\frac{1}{b}\sum_{i=1}^b
\frac{ \sum_{\pi\in
\Pi} \sum_{m=1}^{l}a_{i\pi(j_m)} \exp(t^T a_{i\pi})}{\sum_{\pi
\in\Pi
}\exp(t^T a_{i\pi})}<\sum_{j=k-l+1}^{k}
\bar{A}_{j}.
\end{equation}
So $\Omega\subset\mathcal{P}$.

Let us prove that $\Omega$ is a convex set. Let $x, y \in\Omega$ and
$c + d =1$, $c$, $d > 0$. Then for all $ t \in\mathbf{R}^{k-1}$
\[
t^T (c x + d y) -\ka(t) = c \bigl(t^T x-\ka(t)\bigr) + d
\bigl(t^Ty -\ka(t)\bigr) \le c \La(x) + d \La(y) < \infty.
\]
Since the expression $t^T (c x + d y) -\ka(t)$ is bounded and convex,
it has a maximum, so that $c x + d y \in\Omega$ and
\[
\La(c x + d y) \le c \La(x) + d \La(y),
\]
so $\Omega$ is convex.

To see that each vertex of the polytope is a limiting point of $\Omega
$, consider any vertex $\bar{A}_{\hat{\pi}} = (\bar{A}_{\hat{\pi}(1)},
\bar{A}_{\hat{\pi}(2)}, \ldots, \bar{A}_{\hat{\pi}(k-1)})$.
Suppose $\hat{\pi}(k)= j$ 
and define $l_1,\ldots,l_k$ such that $\hat\pi(l_i)=i$, so $l_j=k$.
We can show that
%
\begin{equation}
\label{limmit} \lim_{t_{l_{j+1}} \to\infty} \cdots\lim_{t_{l_k} \to\infty}
\lim_{t_{l_{j-1}} \to-\infty} \cdots\lim_{t_{l_1} \to-\infty} \ka'
(t) = (\bar{A}_{\hat{\pi}(1)}, \ldots, \bar{A}_{\hat{\pi}(k-1)}).
\end{equation}
To see this, write
%
\begin{eqnarray}
\label{limfin}
\ka' (t) &=& \frac{1}{b}\sum
_{i=1}^b \frac{ \sum_{\pi\in\Pi}
a_{i\pi}
\exp(t^T (a_{i\pi}-a_{i\hat{\pi}}))}{\sum_{\pi\in\Pi}\exp(t^T
(a_{i\pi}-a_{i\hat{\pi}}))}
\nonumber
\\[-8pt]
\\[-8pt]
\nonumber
&= &\frac{1}{b}\sum_{i=1}^b
\frac{ a_{i\hat{\pi}}+\sum_{\pi\in
\Pi, \pi
\neq \hat{\pi}} a_{i\pi} \exp(t^T (a_{i\pi}-a_{i\hat{\pi
}}))}{1+\sum_{\pi\in\Pi, \pi\neq \hat{\pi}}\exp(t^T (a_{i\pi}-a_{i\hat
{\pi}}))}.
\end{eqnarray}
Note that $a_{i\hat{\pi}(l_1)}=a_{i1}$ is the smallest entry in the
$i$th row so the coefficient of $t_{l_1}$, $a_{i\pi(l_1)}-a_{i\hat
{\pi
}(l_1)} = a_{i\pi(l_1)}-a_{i1}$, is either positive or zero for any
$\pi\in\Pi$. As $t_{l_1} \to-\infty$, only the permutations with
$\pi
(l_1) = 1$ give nonzero terms, so
%
\begin{equation}\qquad
\lim_{t_{l_1} \to-\infty} \ka' (t) = \frac{1}{b}\sum
_{i=1}^b \frac{
a_{i\hat{\pi}} + \sum_{\pi\in\{\Pi\dvtx \pi(l_1) = 1\}, \pi\neq
\hat{\pi
}} a_{i\pi} \exp(t^T (a_{i\pi}-a_{i\hat{\pi}}))}{1+\sum_{\pi\in
\{\Pi
\dvtx \pi(l_1) = 1\}, \pi\neq \hat{\pi}}\exp(t^T (a_{i\pi}-a_{i\hat
{\pi
}}))}.\label{lim1}
\end{equation}
Continuing to take limits, in the order given in (\ref{limmit}),
removes all but the first term in the numerator and denominator of
(\ref
{lim1}), to prove (\ref{limmit}).
Thus, the vertex $\bar{A}_{\hat{\pi}}$ is a limiting point of
$\Omega$.
Since $\hat\pi$ is arbitrary, this holds for each vertex of the
polytope. Since $\Omega$ is a convex set enclosed by the edges of the
polytope, $\Omega$ is the interior of the polytope.
\end{pf*}

\begin{pf*}{Proof of Theorem~\ref{Th2}}
Consider the expression $t^Tx - \ka(t)$. Using previous notation set
$x=\bar{A}_{\hat{\pi}}$. Using the definition (\ref{emcum}), we can write
%
\begin{equation}
\label{vertexlim} t^T\bar{A}_{\hat{\pi}} - \ka(t) = -
\frac{1}{b}\sum_{i=1}^b \log
\frac{1}{k!} \sum_{\pi\in\Pi}e^{t^T
(a_{i\pi}-a_{i\hat{\pi}})},
\end{equation}
since $\bar{A}_{\hat{\pi}} = \frac{1}{b}\sum_{i=1}^b a_{i\hat{\pi}}$.
Then by the same argument used in Theorem~\ref{theo2}, we get
\[
\lim_{t_{l_{j+1}} \to\infty} \cdots\lim_{t_{l_k} \to\infty} \lim
_{t_{l_{j-1}} \to-\infty} \cdots\lim_{t_{l_1} \to-\infty} \bigl[ t^T
\bar{A}_{\hat{\pi}} - \ka(t)\bigr]= \log k!.
\]
From the definition of supremum and equation (\ref{emlam}),
\[
\La(x) = \sup_{t} \bigl\{t^T
\bar{A}_{\hat{\pi}} - \ka(t)\bigr\} = \log k!.
\]
Since $\hat{\pi}$ is chosen arbitrarily $\La(x)$ is equal to $\log k!$
on all vertexes. These are the extreme points of $\Omega$ and $\Lambda
(x)$ is convex, so $\Lambda(x)$ takes its maximum on the vertexes and
is finite on all points of $\mathcal{P}$ and so on the boundary of
$\Omega$.
\end{pf*}

\begin{pf*}{Proof of Theorem~\ref{Th3}} Using (\ref{bdyeqs}), choose $x$ so that
\[
\sum_{j=1}^l x_{s_j} = \sum
_{j=1}^l\bar{A}_{j}
\]
is true for some $\{s_j\}$ and $l$. The alternative choice will follow
in the same way. So
\[
x^T t =\sum_{j=1}^{l-1}
x_{s_j} (t_{s_j} -t_{s_l})+ \sum
_{j=1}^l\bar{A}_{j} t_{s_l}+
\sum_{j=l+1}^{k-1} x_{s_j}
t_{s_j}
\]
and, from (\ref{emcum}), $\ka(t)$ can be written
\[
\frac{1}{b}\sum_{i=1}^b \log
\frac{1}{k!} \sum_{\pi\in\Pi} \exp\Biggl({\sum
_{j=1}^{l-1} a_{i\pi(j)} (t_{s_j}-t_{s_l})
+\sum_{j=1}^{l} a_{i\pi(j)}
t_{s_l}+ \sum_{j=l+1}^{k-1}
a_{i\pi(j)} t_{s_j}}\Biggr).
\]
Make the substitution
\[
u_j = \cases{ %
t_{s_j} -
t_{s_l},&\quad $\mbox{for } 1\le j < l,$
\vspace*{2pt}\cr
t_{s_j},&\quad $\mbox{for } l \le j\le k-1,$}
\]
and use the first equality in (\ref{bdyeqs}), to write $x^T t -\ka
(t)$ as
\[
\sum_{j=1,j\neq l}^{k-1} x_{s_j}
u_j-\frac{1}{b}\sum_{i=1}^b
\log \frac
{1}{k!} \sum_{\pi\in\Pi} \exp \Biggl(\sum
_{j=1,j\neq l}^{k-1} a_{i\pi(j)}
u_j+u_l\sum_{j=1}^{l}
(a_{i\pi(j)}-a_{ij}) \Biggr),
\]
where for each $1\le j \le l$, $\sum_{j=1}^l (a_{i\pi(j)}-a_{ij}) \ge
0$. Let $u_l \to-\infty$ and we have
\[
\lim_{u_l \to-\infty}\bigl(x^Tt - \ka(t)\bigr) = \sum
_{j=1,j\neq l}^{k-1} x_{s_j}
u_j-\frac{1}{b}\sum_{i=1}^b
\log\frac{1}{k!} \sum_{\pi\in
\hat
{\Pi}} \exp \Biggl(\sum
_{j=1,j\neq l}^{k-1} a_{i\pi(j)} u_j
\Biggr),
\]
where $\hat{\Pi} = \{ \pi\in\Pi\dvtx  \sum_{j=1}^l (a_{i\pi
(j)}-a_{ij}) =
0\}$. Let $\hat{\Pi}_1$ and $\hat{\Pi}_2$ be sets of all permutations
$\hat{\pi}_1$ and $\hat{\pi}_2$ of integers $\{ 1, \ldots, l\}$ and
$\{
l+1, \ldots, k\}$, respectively. Then the above expression can be rewritten
%
\begin{eqnarray}
\label{tt} \lim_{u_l \to-\infty}\bigl(x^Tt - \ka(t)\bigr)
&= &\sum_{j=1}^{l-1} x_{s_j}
u_j-\frac{1}{b}\sum_{i=1}^b
\log\frac
{1}{l!} \sum_{\hat{\pi}_1 \in\hat{\Pi}_1} e^{\sum_{j=1}^{l-1} a_{i\hat{\pi}_1(j)} u_j}
\nonumber
\\
&&{}+ \sum_{j=l+1}^{k-1} x_{s_j}
u_j -\frac{1}{b}\sum_{i=1}^b
\log \frac
{1}{(k-l)!}\sum_{\hat{\pi}_2 \in\hat{\Pi}_2} e^{\sum_{j=l+1}^{k-1}
a_{i\hat{\pi}_2(j)} u_j}
\\
&&{} +\log\frac{k!}{l!(k-l)!}.\nonumber
\end{eqnarray}
Now taking suprema over $u_1,\ldots,u_{l-1}$ and $u_{l+1},\ldots
,u_{k-1}$ in the first two terms on the right in (\ref{tt}), the
statement of the theorem follows.
\end{pf*}

%



\printaddresses
\end{document}